\documentclass[12pt]{amsart}

\usepackage{amsmath,amssymb,amscd,amsxtra,enumerate}
\usepackage{graphics}
\usepackage{latexsym}
\usepackage{mathtools}
\usepackage{cases}
\usepackage{comment}
\usepackage{ulem}
\usepackage[pdftex]{graphicx,color}

\usepackage{caption}
\usepackage{csquotes}

\newtheorem{thm}{Theorem}
\newtheorem{lemma}{Lemma}
\newtheorem{prop}{Proposition}

\theoremstyle{remark}
\newtheorem{rem}{Remark}
\theoremstyle{definition}

\usepackage[margin=1in]{geometry}

\newcommand{\R}{\mathbb{R}}

\newcommand{\Z}{\mathbb{Z}}

\newcommand{\cG}{\mathcal{G}}

\newcommand{\ip}[2]{\langle#1,#2\rangle}

\newcommand{\e}{\mathrm{e}}

\usepackage[numbered,framed]{matlab-prettifier}
\begin{document}

\title{On the frame set of the second-order B-spline}

\author{A.~Ganiou D.~Atindehou}
\address{ Département de Math\'ematiques/FAST/UAC, 01 BP : 4521, Cotonou 01, B\'enin}
\email{ganiou.atindehou@fast.uac.bj}

\author{Christina~Frederick}
\address{ Department of Mathematical Sciences, New Jersey Institute of Technology}
\email{christin@njit.edu}

\author{Y\'eb\'eni B.~Kouagou}
%u\address{ Institut de Math\'ematiques et de Sciences Physiques (IMSP), 01 BP 613, Porto-Novo, B\'enin} \email{yebeni.kouagou@gmail.com}

\author{Kasso A.~Okoudjou}
\address{Department of Mathematics, Tufts University, Medford, MA 02155}
\email{kasso.okoudjou@tufts.edu}
\thanks{Y\'eb\'eni B.~Kouagou suddenly passed away in 2018, a few weeks after the first version of this work was released. This version is dedicated to his memory.}


\date{\today}


\begin{abstract}
The frame set of a function  $g\in L^2(\R)$ is the set of all parameters $(a, b)\in \R^2_+$ for which the collection of time-frequency shifts of $g$ along $a\Z\times b\Z$ form a Gabor frame for $L^2(\R).$ Finding the frame set of a given function remains a challenging open problem in time-frequency analysis. In this paper, we establish new regions  of the frame set of the second-order B-spline. Our method uses the compact support of this function to partition a subset of the putative frame set and finds an explicit dual window function in each subregion. Numerical evidence indicates the existence of further regions belonging to the frame set.

\end{abstract}

\maketitle \pagestyle{myheadings} \thispagestyle{plain}
\markboth{A. G. D.  Atindehou, C. Frederick, Y. B. Kouagou, And K. A. Okoudjou}{On The Frame Set Of The Second-Order B-Spline}

%\tableofcontents

\section{Introduction and main results}\label{sec1}
Given a window $g \in L^2(\R)$, and  $a, b>0$, the collection of time-frequency shifts 
\begin{equation*}
  \cG(g,a, b)=\left\{M_{\ell b}T_{ka}g =\e^{2\pi i\ell b\cdot}g(\cdot-ka):(\ell,k)\in \mathbb{Z}^2\right\}
\end{equation*} is called a Gabor frame for $L^2(\R)$ if there exist constants $A, B>0$ such that 
\begin{equation*}
 \label{frame inequality}
 A \|f\|^2_2\leq\sum_{\ell,k\in \Z}  |\ip{f}{M_{\ell b}T_{ka}g}|^2\leq B\|f\|^2_2,
 \end{equation*}
 for all $f\in L^2(\R)$. When $\cG(g,a, b) $ is a Gabor frame, there exists a dual window $\gamma \in L^2(\R)$ such that $\cG(\gamma,a, b)$ is also a Gabor frame for $L^2(\R)$ called the (canonical) dual to $\cG(g, a, b)$. Consequently,  for any $f\in L^2(\R)$ we have the following reconstruction formulas:
 
 $$f=\sum_{k \ell \in \Z}\ip{f}{M_{\ell b}T_{k a}\gamma} M_{\ell b}T_{k a}g =\sum_{k \ell \in \Z}\ip{f}{M_{\ell b}T_{k a} g} M_{\ell b}T_{k a}\gamma.  $$ We refer to  \cite{Ole4, Groc2001} for more on Gabor frame theory.

Despite the outstanding advances in the theory and applications of Gabor frames over the last three decades, the problem of characterizing the set of all points $(a, b)\in \R^2_+$ such that $\cG(g, a, b)$ is a frame for a given $g\in L^2(\R)$ remains largely unresolved.  This set is  known as the {\it frame set} of $g$ and will be denoted by $\mathcal{F}(g)$.  The current state-of-the-art result in this direction states that if $g$ is either:
\begin{enumerate}
    \item in $\{e^{-\pi x^2}, \frac{1}{\cosh x}, \chi_{[0, \infty)}e^{-x}, e^{-|x|}\},$ or
    \item a totally positive function of finite type, or
    \item a totally positive function of exponential type 
\end{enumerate}
then $ \mathcal{F}(g)=\{(a, b)\in \R^2_+, ab<1\}$ (see \cite{Groch, Grosto, JanStroh, Jan, Jans, Janss, Lyuba, Seip, SeiWal}). At the same time,  the frame set of $g=\chi_{[0, c]}$, known as the Janssen's tie \cite{Jans}  has a more complex structure that was fully described by Dai and Sun \cite{DaiSun, GuHan}. 

\begin{figure}
 \includegraphics[width=.8\textwidth]{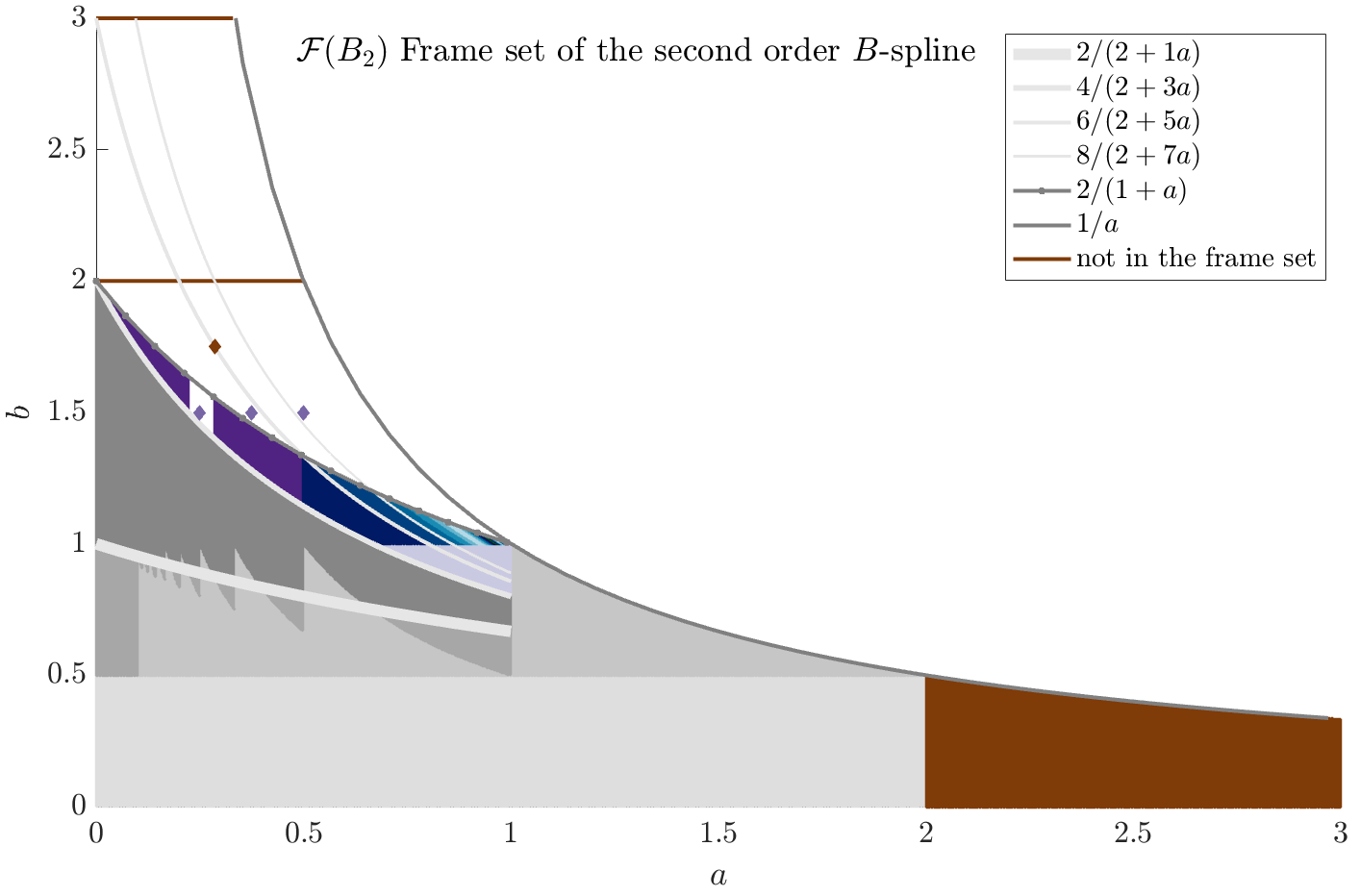}
 \caption{ Sketch depicting known results on the frame set of $B_2$. The shaded regions in gray were proved to be in $\mathcal{F}(B_2)$ in \cite{AtiKouOko1,Olehon, Olehon1, Ole4, Grojan, Lemniel, Kloosto} and the light purple region was established as a subset of $\mathcal{F}(B_2)$ in \cite{AtiKouOko1}. Points that are brown are known not to be in $\mathcal{F}(B_2)$. Theorem \ref{thm-meta} asserts that the union of $\Gamma_3$ (dark purple) and $\cup_{m=3}^{\infty}\Lambda_m$ (shades of blue) belong to $\mathcal{F}(B_2)$. We show that the additional points marked with purple diamonds are in $\mathcal{F}(B_2)$.}
\label{fig:figure1}
\end{figure}

In this paper we consider the frame set of the 
 $B$-spline of order $2$: 
 $$B_2(x)=\begin{cases}
  1+x &  x\in[-1, 0] \\
  1-x & x\in[0,1].\end{cases}$$ 
 It is known that $\mathcal{F}(B_2)$ is an open set in $\R^2_+$ \cite{FeiKai04}, but a full characterization of this set remains an open question. To date, it has been shown that   $$\left\{(a,b)\in\mathbb{R}^2_+:ab<1, 0<a<2,0<b\leq \max_{a}\left(1, \frac{4}{2+3a}\right)\right\} \subset \mathcal{F}(B_2)$$ (see Figure \ref{fig:figure1} for a sketch and \cite{AtiKouOko1,Olehon, Olehon1, Ole4, Grojan, Lemniel, Kloosto}  for details).

 In an earlier work \cite{AtiKouOko1}, we introduced a framework for determining the frame sets of compactly supported functions, including the B-splines of order $N\geq 2$. This framework unified many of the known results on frame sets of $B$-splines \cite{Olehon, Olehon1, Ole4, Grojan, Lemniel, Kloosto}. 
In this paper, we use a similar linear algebra based approach to shed new light on set $\mathcal{F}(B_2)$ and prove the following result

\begin{thm}\label{thm-meta} The frame set of the second-order $B$-spline  contains the set $\Gamma_3\cup \Lambda$, that is,
 $$\mathcal{F}(B_2) \supset \Gamma_3 \cup \Lambda, \qquad \Lambda:=  \bigcup^{\infty}_{m=3}{\Lambda_m},$$
where \begin{align*}
  \Gamma_3&:=\left\{(a,b)\in\mathbb{R}^2_+:a\in \left(0, \frac{2}{9}\right]\cup \left(\frac{2}{7}, \frac{1}{2}\right), b\in \left(\frac{4}{2+3a}, \frac{2}{1+a}\right]\right\}
  \\
  \Lambda_3&:=\left\{(a,b)\in\mathbb{R}^2_+:a\in \left[\frac{1}{2}, \frac{4}{5}\right], b\in \left(\frac{4}{2+3a}, \frac{6}{2+5a}\right], b>1\right\},
  \end{align*} and for $m\geq 4$
 
 \begin{align*}\label{setlbdam}
  \Lambda_m=\left\{(a,b)\in\mathbb{R}^2_+\right.:a&\in \left[ \frac{m-3}{m-2}, \frac{2(m-1)}{2m-1}\right]\\
  b&\left.\in\left( \frac{2(m-1)}{2+(2m-3)a},
   \min_{a}\left( \frac{2m}{2+(2m-1)a}, \frac{2}{1+a}\right)\right], b>1\right\}.
  \end{align*}  
 \end{thm}

The sets $\Lambda_m$, $m\geq 3$ and $\Gamma_3$ appear naturally in the region $\{(a, b) \in \R^2_+: ab<1\}$ due to the compact support of the second-order $B$-spline. Figure~\ref{fig:figure1} displays these sets in varying shades of blue. In particular, the dark purple and blue regions in Figure \ref{fig:figure1} illustrates Theorem~\ref{thm-meta}. While the methods presented here do not fully characterize $\mathcal{F}(B_2)$, they enable us to numerically find new putative points in this set.

The proof of Theorem~\ref{thm-meta}  is divided in two parts: Theorem~\ref{thm1:b2new} (concerning the regions $\Lambda_m$) and Theorem \ref{thm2:b2new} (concerning the region $\Gamma_3$), and is based  on methods  developed in 
\cite[Proposition 2]{AtiKouOko1} which give a necessary and sufficient condition on two Bessel Gabor systems to be dual. We outline the method before giving its full details in the next two sections. 

To deal with the regions $\Lambda_m$, $m \geq 3$,%,  $a\in (0,2)$, and  $b\in (\frac{2(m-1)}{2+(2m-3)a}<\frac{2m}{2+(2m-1)a}]$.
 let $h$ be a bounded, real-valued function with support on 
 $[-\frac{(2m-1)a}{2},\frac{(2m-1)a}{2}]$.
 Then, it is known that the Gabor systems $\cG(B_2,a, b)$ and $\cG(h,a, b)$ are dual frames for $L^2(\mathbb{R})$ if and only if 
 \begin{equation}\label{eq:dualed}
  \sum_{k=1-m}^{m-1}{B_2(x-\ell/b+ka)h(x+ka)}=b\delta_{\ell,0},\qquad |\ell|\leq m-1\, , \  \mbox{for a.e}\ x\in[-\frac{a}{2},\frac{a}{2}].
  \end{equation}
Our main goal is to show the existence of a bounded compactly function $h$ that solves~\eqref{eq:dualed} when $(a, b)\in \Lambda_m$, for $m\geq 3$. 

To do this, we rewrite \eqref{eq:dualed} as a matrix-vector equation using the compact support of $B_2$,
\begin{equation}\label{eq-matrix}
 G_mH_m=
\begin{pmatrix}F_{m-1} &A_{m-1}\\0 & C_{m-1}\end{pmatrix}H_m=V
\end{equation}
where $G_m$ is a $(2m-1)\times (2m-1)$ matrix, and the column vectors $H_m$ and $V$ are given by \[H_m(x)=\left[h(x+ka)\right]_{|k|\leq m-1}, \qquad V=[b\delta_{\ell,0}]_{|\ell|\leq m-1}.\] Furthermore, 
$C_{m-1}$ is 
a $(m-2)\times(m-2)$ upper triangular matrix with strictly positive diagonal entries, and  $F_{m-1}$ is the $(m+1)\times(m+1)$ tridiagonal matrix
given by 
\begin{equation}
 F_{m-1}(x)=\left[B_2(x-\frac{\ell}{b}+ka)\right]_{1-m\leq\ell,k\leq 1}.\label{eq:Fm}
\end{equation} 
We observe that  the diagonal entries of $F_{m-1}$ are also strictly positive. Therefore, the existence of a dual is guaranteed when $G_m$ has a strictly positive determinant, or equivalently (because  of its  structure), 
that the tridiagonal matrix $F_{m-1}$ has a strictly positive determinant.

Although recursive formulas for computing the determinants of tridiagonal matrices exist, the present case requires more than showing these determinants are nonzero. As such, we exploit the properties of $B_2$ to prove recursively that  the determinant of $F_{m-1}$ is strictly positive. This allows us to solve explicitly for the bounded and compactly supported function $h$ that generates a Bessel system and is dual to $g$.  We carry out this analysis in  Section~\ref{sec2}, establishing the invertibility of the matrices when $(a, b) \in \Lambda$. 

To deal with the region $\Gamma_3$, we effectively compute the determinant a  $4\times 4$ matrix that models ~\eqref{eq:dualed}. However, in this case the matrix is no longer tridiagonal. Nonetheless, we again exploit the structure of $B_2$ to show that the determinant of this matrix is also strictly positive which allows us to find a compactly supported and bounded function $h$ which generates a dual Gabor frame to $\cG(B_2,a, b)$. This is established in 
 Section~\ref{subsec2.2}. 

{\bf Notation.} In the sequel, given a $p\times p$ matrix $A$ and $E\subset \{1, 2, 3, \hdots, p\}$, we denote by $A^{E}$ the $\# E \times \# E$  sub-matrix of $A$  using rows and columns from $E$, and denote by $|A|$ the determinant of the matrix $A$. We also denote by $g_{l,k}$ the function
\[g_{l,k}(x):= B_2(x-\frac{l}{b}+ka).\]

\section{$\mathcal{F}(B_2)$ contains $\Lambda$}\label{sec2}
In this section we prove the first part of Theorem~\ref{thm-meta} by establishing the following result. 

\begin{thm}\label{thm1:b2new}
 For $m\geq3$, let $(a,b)\in \Lambda_m$. Then, the Gabor system $\cG(B_2,a, b)$ is a  frame for $L^2(\mathbb{R})$, and there is
a unique dual window $h \in L^2(\mathbb{R})$ such that supp$(h)\subseteq\left[-\frac{2m-1}{2}a,\frac{2m-1}{2}a\right]$.
Furthermore, for each $(a,b)\in \Lambda$, the Gabor system $\cG(B_2,a, b)$ is a  frame for $L^2(\mathbb{R})$.
\end{thm}

To prove Theorem~\ref{thm1:b2new} we only need to show that ~\eqref{eq:dualed} has a solution $h$ that is a bounded and  compactly supported function. As mentioned earlier, the determinant of the (block) matrix $G_{m}$ is $$ |G_{m} |= |F_{m-1}| |C_{m-1} |= |F_{m-1} | \prod_{1-m}^{-2}{g_{-k,-k}}$$ where we used the fact $C_{m-1}$ is an upper triangular matrix. Because for each $k=1-m, \hdots, -2$, $g_{-k,-k}>0$ on $[-\frac{a}{2}, 0]$,  we turn our attention to establishing that the determinant  $ |F_{m-1}|$ of this tridiagonal matrix is strictly positive when  $(a,b)\in\Lambda_m$.

\begin{prop} \label{prop:determinantFm} For $m\geq3$, let $(a,b)\in\cup_{k=m}^{\infty}\Lambda_k$. Then  $ |F_{m-1}|>0$ on $\left[- \frac{a}{2};0\right]$.
\end{prop}

We will prove the result by showing that $|F_{m-1}|$ never vanishes on $ [-\frac{a}{2}, 0]$. Since this matrix is a tri-diagonal, we could rely on standard formulas to find its determinant.  However, the challenging part is to establish that the determinant is never $0$. We will do this by an induction argument on $m$, relying on the fact that $B_2$ is a compactly supported piecewise linear function. In Lemma~\ref{lem:basecase} we prove the result for the base case $m=3$, and complete the induction argument  in Lemma~\ref{lem:inductionpart}.

\begin{figure}
    \centering
    \includegraphics[width=.5\textwidth]{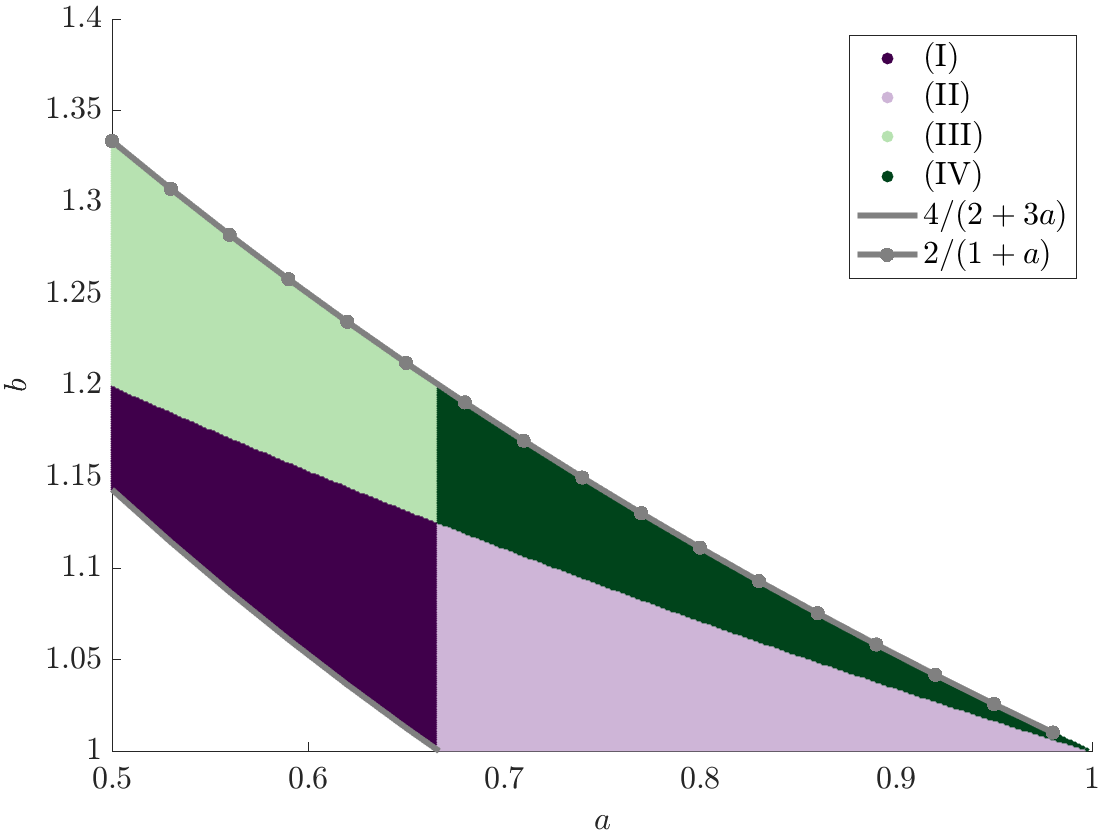}\includegraphics[width=.5\textwidth]{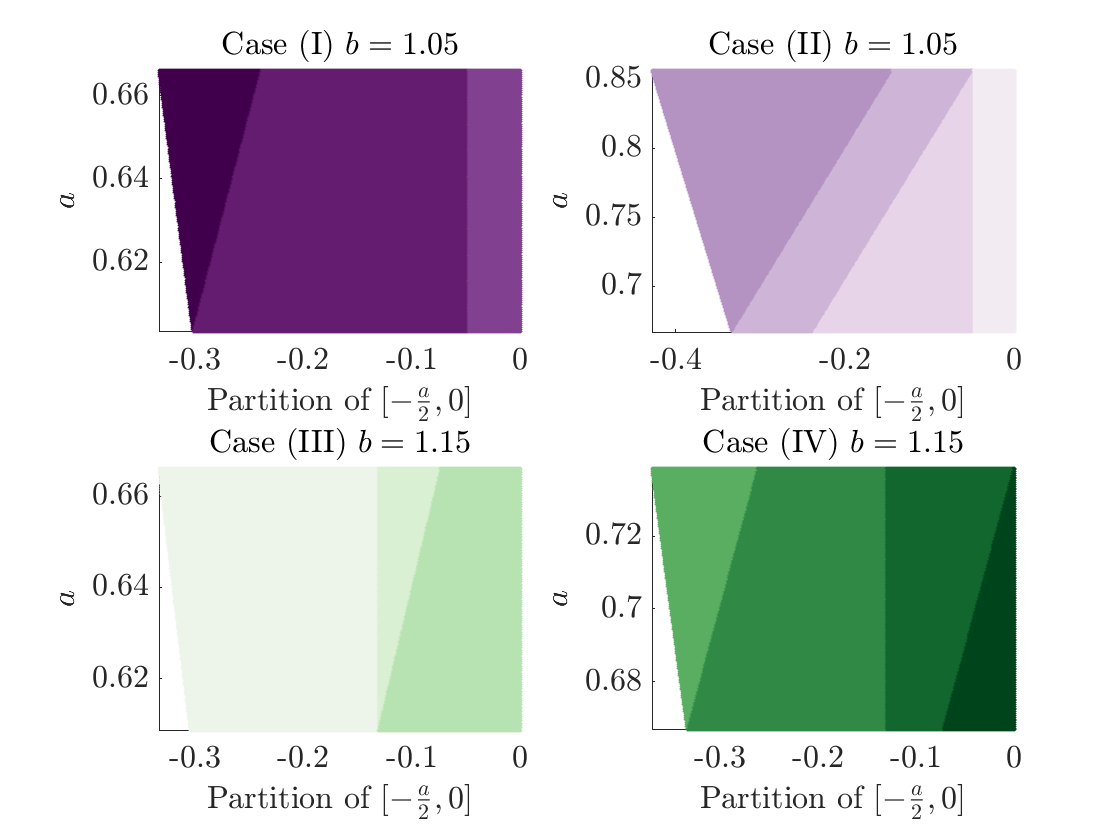}
    \caption{The frame set of $B_2$ established in Lemma \ref{lem:basecase} by partitioning $\cup_{k=3}^{\infty}\Lambda_k$ into the four regions shown in the left plot. In each region,  for each fixed $b$, the interval $[-\frac{a}{2},0]$ is subdivided into three intervals (Cases (I) and (III)) or four intervals (Cases (II) and (IV)) as depicted by the shading in the plots on the right.}
    \label{fig:basecase}
\end{figure}
\begin{lemma}\label{lem:basecase}  Let  $(a, b) \in \cup_{k=3}^\infty\Lambda_k.$ %=\{(a, b):a\in [\frac{1}{2},1],\, b\in (\frac{4}{2+3a},  \frac{2}{1+a}]\,\, ,{\text and}\,\,  b>1\}.$ 
Then, $|F_{2}|>0$ on $[-\frac{a}{2}, 0]$.

\end{lemma}

\begin{proof} We note that for $ (a, b) \in \cup_{k=3}^\infty\Lambda_k$

\begin{equation}\label{matrixB2} 
F_2=\begin{pmatrix}
  g_{-2,-2} &g_{-2,-1}&0 & 0\\
 g_{-1,-2} &g_{-1,-1}&g_{-1,0}&0\\ 
  0&g_{0,-1} &g_{0,0}&g_{0,1}\\
  0 & 0 &g_{1,0} &g_{1,1}\\
 \end{pmatrix}.
 \end{equation}
In each case below, we write $[- \frac{a}{2},0]$ as a union of intervals and show that $|F_2|>0$ on each of them. Figure \ref{fig:basecase} shows plots of the four regions of $\cup_{k=3}^\infty\Lambda_k$ considered (left) as well as the partitions of the interval $[-\frac{a}{2},0]$ (right). %We will denote by the point \[q_{k,l} = 1-\frac{k}{b}+la\]
 
\begin{enumerate}[(I)]
    \item For $a\in [\frac{1}{2}, \frac{2}{3})$ and $2- \frac{3}{b}+a\leq 0$, 
    
    \begin{equation*}
 { [- \frac{a}{2},0]=[- \frac{a}{2},1- \frac{2}{b}+a)\cup 
 [1- \frac{2}{b}+a,-1+ \frac{1}{b}]\cup (-1+ \frac{1}{b},0]}. 
 \end{equation*} 
% \begin{equation*}
%  [- \frac{a}{2},0]=[- \frac{a}{2},q_{2,1})\cup 
%  [q_{2,1},-q_{1,0}]\cup (-q_{1,0},0] 
%  \end{equation*}
 On the first interval, $g_{1,0}= 0$, and since $ g_{-1,-1}>g_{-1,-2}$, $ g_{-1,-1}\geq g_{-1,0}$ and
\[\begin{vmatrix}
 g_{-2,-2} &   0  \\
 0 &g_{0,0} \\
 \end{vmatrix}>\begin{vmatrix}
  g_{-2,-2} &  g_{-2,-1}  \\
  0 &g_{0,-1} \\
 \end{vmatrix}+\begin{vmatrix}
 g_{-2,-1} & 0\\
g_{0,-1} &g_{0,0}\\
 \end{vmatrix},\] it follows that $|F_{2} |=g_{1,1} |F^{\{1,2,3\}}_{2} |>0$. On the second interval, $g_{2,1}=
g_{1,0}=0$, and since $g_{-1,-1}>g_{-1,0}$ and $g_{0,-1}<g_{0,0}$, it follows that
$ |F_{2}|=g_{-2,-2} g_{1,1} |F^{\{2,3\}}_{2} |>0$. On the third interval, $g_{-2,-1}=0$, and since $g_{0,0}>\max\{g_{0,-1},g_{0,1}\}$, and
\[\begin{vmatrix}
 g_{-1,-1} &   0  \\
 0 & g_{1,1} \\
 \end{vmatrix}>\begin{vmatrix}
  g_{-1,-1} &  g_{-1,0}  \\
  0 & g_{1,0} \\
 \end{vmatrix}+\begin{vmatrix}
 g_{-1,0} & 0\\
 g_{1,0} & g_{1,1}\\
 \end{vmatrix},\] it follows that $|F_{2} |=g_{-2,-2} |F^{\{2,3,4\}}_{2} |>0$.

    \item For $a\in [\frac{2}{3}, 1]$ and $2- \frac{3}{b}+a\leq 0$,
\[[-\frac{a}{2}, 0]= [-\frac{a}{2}, a-1]\cup (a-1, 1-\frac{2}{b}+a)\cup [1-\frac{2}{b}+a, \frac{1}{b} -1]\cup (\frac{1}{b}-1, 0].\] 
% \[[-\frac{a}{2}, 0]= [-\frac{a}{2}, -q_{0,-1}]\cup (-q_{0,-1}, q_{2,1})\cup [q_{2,1}, -q_{1,0}]\cup (-q_{1,0}, 0].\] 
On the first interval, we have $g_{0,-1}=g_{1,0}=0$, and since $g_{-2,-2}>g_{-2,-1}$, $g_{-1,-2}< g_{-1,-1}$, it follows that $ |F_{2} |=\prod_{\ell=0}^1 g_{\ell,\ell}  |F_{2}^{\{1,2\}} |>0$. 
The last three intervals are treated as in the previous case.

\item For $ [\frac{1}{2}, \frac{2}{3})$ and $2- \frac{3}{b}+a>0$,
\[ [- \frac{a}{2},0]=[- \frac{a}{2},-1+ \frac{1}{b}]\cup 
(-1+ \frac{1}{b},1- \frac{2}{b}+a)\cup [1- \frac{2}{b}+a,0].\]
%  \[ [- \frac{a}{2},0]=[- \frac{a}{2},-q_{1,0}]\cup 
% (-q_{1,0},q_{2,1})\cup [q_{2,1},0].\]
On the first interval, $g_{1,0}=0$. As in (I), we conclude that  
$ |F_{2} |=g_{1,1} |F^{\{1,2,3\}}_{2} |>0$. On the second interval, we have $g_{0,0}>g_{0,-1}$, $g_{0,0}>g_{0,1}$, $|F^{\{1,2\}}_{2}|>0$ and $|F^{\{1,2,4\}}_{2}|>g_{1,0}|F^{\{1,2\}}_{2}|+
g_{-2,-2}g_{-1,0}g_{1,1}$, and therefore \begin{equation*} 
|F_{2}|=
\begin{vmatrix}
 g_{-2,-2} & g_{-2,-1} &  0 & 0 \\
 g_{-1,-2} & g_{-1,-1} & g_{-1,0} & 0\\
 0 &g_{0,-1} &g_{0,0} &g_{0,1}\\
 0 & 0 & g_{1,0}& g_{1,1}\\ 
 \end{vmatrix}>0.
\end{equation*}
On the third interval, $g_{-2,-1}=0$. As in part (I), we conclude $ |F_{2}|=g_{-2,-2} |F^{\{2,3,4\}}_{2} |>0$. 

    \item For $a\in [\frac{2}{3}, 1]$ and $2- \frac{3}{b}+a>0$, 
 \[
 [- \frac{a}{2},0]=[- \frac{a}{2},-1+a]\cup(-1+a,-1+ \frac{1}{b}]\cup 
 (-1+ \frac{1}{b},1- \frac{2}{b}+a)\cup [1- \frac{2}{b}+a,0].\]
%  $-q_{0,-1} = -1+\frac{0}{b}+a$\\
%  $-q_{1,0} = -1+\frac{1}{b}+0a$\\
%  $q_{2,1} = 1-\frac{2}{b}+a$
%  \[
%  [- \frac{a}{2},0]=[- \frac{a}{2},-q_{0,-1}]\cup(-q_{0,-1},-q_{1,0}]\cup 
%  (-q_{1,0},q_{2,1})\cup [q_{2,1},0].\] 
On the first interval, we have $g_{0,-1}=g_{1,0}=0$, and therefore $|F_{2}|=g_{1,1}g_{0,0}|F^{\{1,2\}}_{2}|>0$. The  remaining three intervals can be treated as above.

\end{enumerate}
\end{proof}

We organize the induction step in the following result.

\begin{lemma}\label{lem:inductionpart} Suppose that for some $m\geq3$, $ |F_{m-1}|>0$ on $[-\frac{a}{2}, 0]$ for each $(a, b)\in \Lambda_k$, for all $k\geq m$.  Then   $ |F_{m} |>0$ on $[-\frac{a}{2}, 0]$ for $(a, b) \in \Lambda_{k},$ and all $ k\geq m+1$.
 
\end{lemma}

\begin{proof} We first prove that  $|F_m|>0$ on $[-\frac{a}{2},0]$ for $(a, b)\in \Lambda_{m+1}$ in the following four cases:

\begin{enumerate}[I)]
\item For $a\in [\frac{m-2}{m-1}, \frac{4m-2}{4m+1}]$ and $2- \frac{3}{b}+a\leq 0$,
   then $-1- \frac{m-2}{b}+(m-1)a\leq- \frac{a}{2}$. In this case, \[\left[- \frac{a}{2},0\right]=
   \left[- \frac{a}{2},1- \frac{m}{b}+(m-1)a\right)\cup\left[1- \frac{m}{b}+(m-1)a,0\right].\]
%   \[\left[- \frac{a}{2},0\right]=
%   [- \frac{a}{2},q_{m, m-1})\cup\left[q_{m, m-1},0\right].\]
On the first interval,  $g_{-k,-(k+1)}=0$ for all $k \in\{-1,...,m-3\}$. Therefore, since 
$g_{-m+1,-m}<  g_{-m+1, -m+1},  
 g_{-m+1, -m+1}> g_{-m+1, -m+2}$,  and 
\begin{align*}
\begin{vmatrix}
 g_{-m,-m}&   0  \\
 0 & g_{-m+2,-m+2} \\
 \end{vmatrix}&>\begin{vmatrix}
 g_{-m,-m+1} &  0  \\
  g_{-m+2,-m+1} & g_{-m+2,-m+2}\\
 \end{vmatrix}+
 \begin{vmatrix}
 g_{-m,-m} & g_{-m,-m+1}\\
 0 & g_{-m+2,-m+1}\\
 \end{vmatrix},\end{align*}
 it follows that
 $|F_m |= \prod_{\ell=-1}^{m-3}{g_{\ell, \ell}} |F^{\{1,2,3\}}_m |>0$. On the second interval, $g_{-m,-(m-1)}=0$. The induction assumption implies that
 $|F_m |=g_{-m,-m}|F_{m-1}|>0$. 

\item For $a\in (\frac{2m-2}{2m-1},  \frac{2m}{2m+1}]$ and $2- \frac{3}{b}+a\leq 0$, we have 
$-\frac{a}{2}\leq-1-\frac{(m-2)}{b}+(m-1)a$, and  
\[ \begin{split} [- \frac{a}{2},0]&=[- \frac{a}{2}, -1- \frac{m-2}{b}+(m-1)a]\\ &\cup(-1-
 \frac{m-2}{b}+(m-1)a,1- \frac{m}{b}+(m-1)a) \\  &\cup [1- \frac{m}{b}+(m-1)a,0].\end{split}\]
%  \[ [- \frac{a}{2},0]=[- \frac{a}{2}, -q_{-(m-2),-(m-1)}]\cup(-q_{-(m-2),-(m-1)},q_{m,m-1})\cup[q_{m,m-1},0]\]
 On the first interval,  $g_{-k,-(k+1)}=0$ for all $k \in\{-1,...,m-2\}$. Since $F^{\{1,2\}}_m$ is diagonally dominant, it follows that
$|F_m |= \prod_{\ell=-1 }^{m-2}{g_{-\ell, -\ell}} |F^{\{1,2\}}_m |>0$. On the second interval,  $g_{-k,-(k+1)}=0$ for all $k \in\{-1,...,m-3\}$. As in part (I), we conclude $
  |F_m |= \prod_{\ell=-1}^{m-3}{g_{-\ell, -\ell}} |F^{\{1,2,3\}}_m |>0.$ On $\left[1- \frac{m}{b}+(m-1)a,0\right]$, $g_{{-m},-(m-1)}=0$. Consequently 
 $ |F_m|=g_{{-m},-m}|F_{m-1} |>0$ by the induction assumption.

\item For $a\in (\frac{4m-2}{4m+1}, \frac{2m-2}{2m-1})$ and $2- \frac{3}{b}+a\leq 0$, the quantity $-1-(m-2)/b+(m-1)a+\frac{a}{2}$ can be either positive or negative, falling into the categories of (I) and (II).

\item If $2- \frac{3}{b}+a>0$, then  $a\in [\frac{m-2}{m-1}, \frac{4m-2}{4m+1}]$ and 
\[\begin{split}
\left[- \frac{a}{2},0\right]&=\left[- \frac{a}{2},-1- \frac{m-3}{b}+(m-2)a\right]\\&\bigcup \left(-1- \frac{m-3}{b}+(m-2)a,1- \frac{m}{b}+(m-1)a\right)\\& \bigcup\left[1- \frac{m}{b}+(m-1)a,0\right].\end{split}\]

On the first interval, $g_{-k,-(k+1)}=0$ for all $k \in\{-1,...,m-3\}$. A in part (I), we have
$|F_m |= \prod_{\ell=-1 }^{m-3}{g_{-\ell, -\ell}} |F^{\{1,2,3\}}_m |>0$. On the second interval, $g_{-k,-(k+1)}=0$ for all $k \in\{-1,...,m-4\}$. We then have $F^{\{1,2,3,4\}}_m(x)=F_2(x+\frac{m-2}{b}-(m-2)a)$, where $F_2$ is the $4\times 4$ matrix~\eqref{matrixB2}. Hence,
\[|F_m |= \prod_{\substack{\ell=-1\\ \ell-m\neq-3,-2,-1,0}}^{m}g_{-\ell, -\ell} |F^{\{1,2,3,4\}}_m |>0.\] On $\left[1- \frac{m}{b}+(m-1)a,0\right]$, $g_{-m,-(m-1)}=0$. Hence,
$ |F_m|=g_{-m,-m} |F_{m-1}|>0$ by the induction assumption. 
\end{enumerate}

 To establish the result for $(a, b)\in \Lambda_k$, $k\geq m+2$, we prove that
  $|F_{m}|>0$ on each interval in a partition of $\left[- \frac{a}{2},0\right]$ and reduce the analysis to the case $(a, b) \in \Lambda_{m+1}$. We omit details of the proof, only indicating the relevant partitions of $[-\frac{a}{2}, 0]$.  
  
 For $2- \frac{3}{b}+a\leq 0$ we have the following partition: $$\left[ \frac{k-3}{k-2}, \frac{2(k-1)}{2k-1}\right]=\left[ \frac{k-3}{k-2}, \frac{4k-6}{4k-3}\right] \bigcup \left( \frac{4k-6}{4k-3}; \frac{2k-4}{2k-3}\right) \bigcup \left[ \frac{2k-4}{2k-3}, \frac{2(k-1)}{2k-1}\right].$$ 
   
   \begin{enumerate}[(I)]  
   \item If $a\in\left[ \frac{k-3}{k-2}, \frac{4k-6}{4k-3}\right]$, then $-1- \frac{k-3}{b}+
   (k-2)a\leq- \frac{a}{2}$, and we write
   \[[- \frac{a}{2},0]= \bigcup_{\ell=1}^{k-2}{(Q_\ell\cup T_\ell)} \cup
   T_{k-1}\] where $T_1=(-1+ \frac{1}{b},0]$,
    $T_{k-1}=[- \frac{a}{2},1- \frac{k-1}{b}+(k-2)a),$  and for $\ell =2, \hdots, k-2$ we have  $
  T_\ell=(-1- \frac{\ell-2}{b}+(\ell-1)a,1- \frac{\ell}{b}+(\ell-1)a);$  while  for $\ell =1, \hdots, k-2$ we have $Q_\ell=[1- \frac{\ell+1}{b}+\ell a,-1- \frac{\ell-2}{b}+(\ell-1)a]$.

 \item If  $a\in\left[ \frac{2k-4}{2k-3}, \frac{2(k-1)}{2k-1}\right]$, then $-1- \frac{k-3}{b}+
 (k-2)a\geq- \frac{a}{2}$, and we have the partition 
 \[[- \frac{a}{2},0]= \bigcup_{\ell=1}^{k-1}{(Q_\ell\cup T_\ell)}\]
 where $ T_1=(-1+ \frac{1}{b},0]$,  for $\ell=2, \hdots, k-1$,  $T_\ell=(-1- \frac{\ell-2}{b}+(\ell-1)a,1- \frac{\ell}{b}+(\ell-1)a)$, while  $\ell=1, \hdots, k-2$,  $Q_\ell=[1- \frac{\ell+1}{b}+ka,-1- \frac{\ell-2}{b}+(\ell-1)a]$, and $Q_{k-1}=[- \frac{a}{2},-1- \frac{k-3}{b}+(k-2)a]$.

 \item  Finally,  if  $a\in\left( \frac{4k-6}{4k-3}; \frac{2k-4}{2k-3}\right)$, then the sign of $-1- \frac{k-3}{b}+(k-2)a+
  \frac{a}{2}$ is not constant. However, by considering the case in which this expression is either positive or negative, we can reduce the analysis to one of the two previous cases. 
\end{enumerate}

For $2- \frac{3}{b}+a>0$, we  have   $a\in \left[ \frac{k-3}{k-2}, \frac{4k-6}{4k-3}\right]$ and
  \[\left[- \frac{a}{2},0\right]=\left( \bigcup_{\ell=1}^{k-1}{\widetilde{Q}_\ell}\right)\cup
  \left( \bigcup_{\ell=2}^{k-1}{\widetilde{T}_\ell}\right)\]
where \[\widetilde{Q}_\ell=\left[1- \frac{\ell+1}{b}+ka,-1- \frac{\ell-3}{b}+(\ell-2)a\right]\text{ and } \widetilde{T}_\ell=(-1- \frac{\ell-3}{b}+(\ell-2)a,1- \frac{\ell}{b}+(\ell-1)a),\] with the convention that
$\widetilde{Q}_{k-1}=[- \frac{a}{2},-1- \frac{k-4}{b}+(k-3)a]$ and $ 
 \widetilde{Q}_1=\left[1- \frac{2}{b}+a,0\right].$

\end{proof}

We can now prove Proposition~\ref{prop:determinantFm}.

\begin{proof}[Proof of Proposition~\ref{prop:determinantFm}]
The result is now proved by induction using Lemma~\ref{lem:basecase} and Lemma~\ref{lem:inductionpart}.
\end{proof}

We are now ready to prove Theorem~\ref{thm1:b2new}.

\begin{figure}
    \centering
    \includegraphics[width=\textwidth]{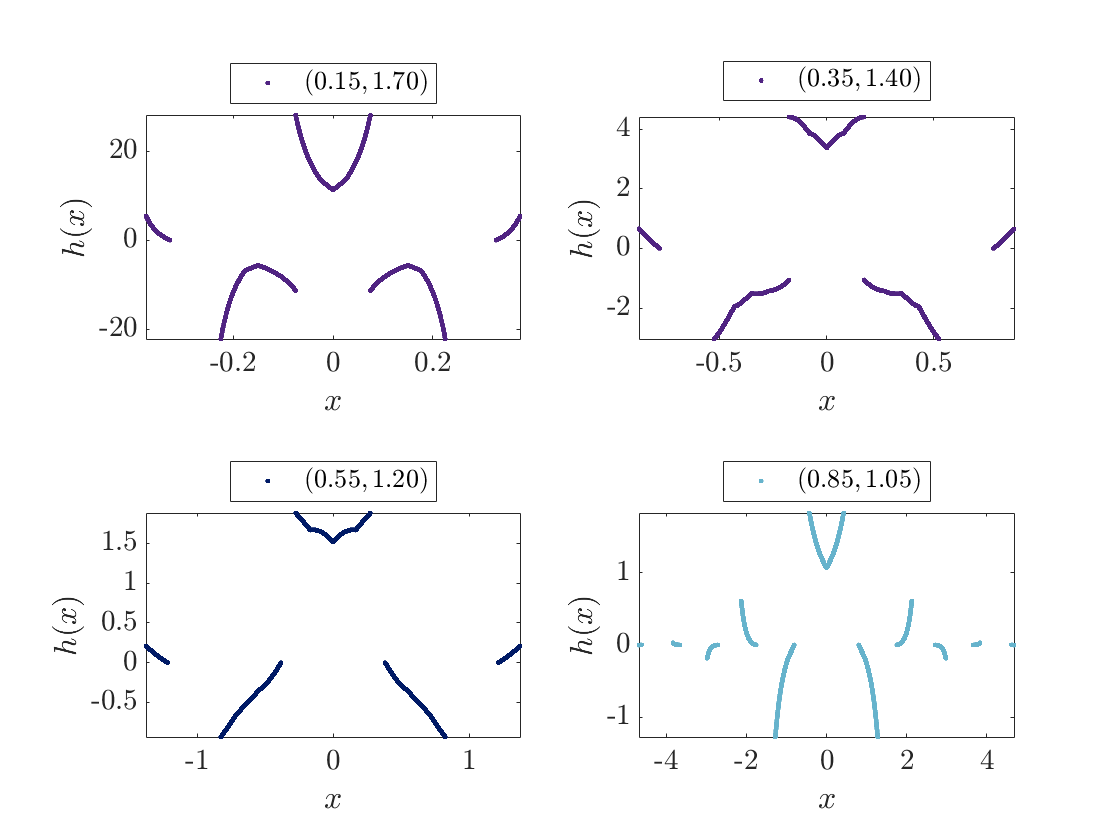}
    \caption{Plots of nonzero values of the dual $h(x)$ of $\mathcal{G}(B_2, a,b)$ for the following $(a,b)$: (0.15, 1.70) and (0.35, 1.40) in  $\Gamma_3$; (0.55, 1.20) in $\Lambda_3$; and (0.85, 1.05) in $\Lambda_6$.}
    \label{fig:duals}

\end{figure}

\begin{proof}[Proof of Theorem~\ref{thm1:b2new}]
 By Proposition~\ref{prop:determinantFm} we know that $G_m$ is invertible. Furthermore,  $|G_m(-x)|=|G_m(x)|$ for all $x\in [-\frac{a}{2}, \frac{a}{2}]$, \cite[Remark 2]{AtiKouOko1}. Let $h$ be defined on $\R$ as follows: let $h(x)=0$ if 
 $x\in\R\setminus\left[- \frac{2m-1}{2}a, \frac{2m-1}{2}a\right]$ and for  $x\in\left[-\frac{2m-1}{2}a,\frac{2m-1}{2}a\right]$ 
let $h(x)$ be the appropriate entry of
$H_m(x)=b(G^{-1}_m(x))_m$, where \begin{equation}H_m(x)=\left[h(x+ka)\right]_{|k|\leq m-1}\label{eq:H}\end{equation} and $(G^{-1}_m(x))_m$ is the $m^{th}$ column vector of the matrix $G^{-1}_m(x).$

Then we can determine $h$ on $(-\frac{a}{2}+ka, ka]$ where $k\in\{1-m,...,m-1\}$. It can be shown that $h$ is even and therefore can be defined on the interval 
 $\left[- \frac{2m-1}{2}a, \frac{2m-1}{2}a\right]$ except at finitely many points (see Figure \ref{fig:duals}).
Since $ |G_m|>0$ on $\left[- \frac{a}{2}, \frac{a}{2}\right]$, we conclude
 that $  \tfrac{1}{|G_m |}$ is a continuous, hence a bounded function on 
 $\left[- \frac{a}{2}; \frac{a}{2}\right]$. Consequently, $h$ is a compactly supported 
and bounded function for which $\cG(h,a, b)$ is a Bessel sequence. By construction, it also follows that $B_2$ and $h$ are dual windows.
\end{proof}

\section{$\mathcal{F}(B_2)$ contains $\Gamma_3$}\label{subsec2.2}
We now turn to the second part of Theorem~\ref{thm-meta} by showing that the following result holds. 

\begin{thm}\label{thm2:b2new}
 Let $(a,b)\in \Gamma_3$. Then, the Gabor system $\cG(B_2,a, b)$ is a  frame for $L^2(\mathbb{R})$, and there is
a unique dual window $h \in L^2(\mathbb{R})$ such that supp$(h)\subseteq\left[-\frac{5a}{2},\frac{5a}{2}\right]$.
\end{thm}

We observe that when $(a,b)\in\Gamma_3$, the matrix $G_3$ becomes
\begin{equation}\label{eq2:keymatrix}
 G_3 =\begin{pmatrix}D  &v \\0 & g_{2,2} \end{pmatrix},
\end{equation}
where $0$ is a $1\times 4$ matrix of $0s$, $v $ is a column vector in $\R^4$ and 
$D $ denotes the $4\times 4$ matrix  obtained by deleting the last row and the last column of $G_3 $ and given by
\begin{equation*}
D=\begin{pmatrix}
g_{-2,-2}  &g_{-2,-1}  & 0 & 0\\
g_{-1,-2}  &g_{-1,-1}  &g_{-1,0}  & g_{-1,1} \\
g_{0,-2}  &g_{0,-1}  &g_{0,0}  &g_{0,1} \\
g_{1,-2}  &g_{1,-1}  &g_{1,0}  &g_{1,1} \\ 
 \end{pmatrix}.
\end{equation*}

Because $g_{2,2} >0$ on $[-\frac{a}{2}, 0]$, 
$G_3 $ is invertible on $[-\frac{a}{2}, 0]$ if and only if $D $ is invertible. 
The following proposition shows that the matrix $D $ is invertible for $(a,b)\in \Gamma_3=\Gamma'_3\cup\Gamma''_3$, where
 \begin{equation}
 \begin{split} 
  \Gamma'_3:=\left\{(a,b)\in\mathbb{R}^2_+:a\in \left(0, \frac{2}{9}\right], b\in \left(\frac{4}{2+3a}, \frac{2}{1+a}\right]\right\}, \text{ and}\\
  \Gamma''_3:=\left\{(a,b)\in\mathbb{R}^2_+:a\in \left(\frac{2}{7}, \frac{1}{2}\right), b\in \left(\frac{4}{2+3a}, \frac{2}{1+a}\right]\right\}. 
  \end{split}
  \end{equation}

\begin{prop}\label{prop:invgamma_3}
Let  $(a, b)\in \Gamma_3=\Gamma'_3\cup\Gamma''_3$. Then $|D |>0$, and therefore $|G_3 |>0$, on  $[-\frac{a}{2}, 0]$.
\end{prop}

\begin{proof}  
We first consider the case $a\in(0, \frac{2}{13}]$. Then $g_{1,-2} >0$ and $-\frac{a}{2}< 1 - \frac{2}{b}+a \leq 0$. A series of computations shows that $|D |>0$, since
\begin{equation}\label{deter1}
 |D (x)|=\begin{cases}
  -\frac{4a}{b^2}(bx-b+2)[(1-b)x+a(1-b)]& x\in [-\frac{a}{2}, 1-\frac{2}{b}+a)\\
   \frac{4a(b-1) (x+a)}{b}  g_{-2,-2}(x)   & x\in [1-\frac{2}{b}+a, 0].\end{cases}
\end{equation}

Next, assume that $a\in(\frac{2}{13}, \frac{1}{5})$. Then $g_{1,-2}(x) >0$ if and only if $x\in (-1+\frac{1}{b}+2a, 1+\frac{1}{b}+2a)$. 

We first treat the subcase $-1+\frac{1}{b}+2a\geq0$, then $g_{1,-2} =0$. Therefore
\begin{equation}\label{deter2}|D(x) |=\begin{cases}
  I(x) & x\in [-\frac{a}{2}, 1-\frac{2}{b}+a)\\
   \frac{4a(b-1) (x+a)}{b}  g_{-2,-2}(x)   & x\in [1-\frac{2}{b}+a, 0].
     \end{cases}
\end{equation}
We now prove that  $I >0$ on the interval $[-\frac{a}{2}, 1-\frac{2}{b}+a)$ by showing that $I' >0$ on  $\left[-\frac{a}{2},0\right]$ and  $I(-\frac{a}{2})>0$. It can be proved that  $f(a,b):=I(-\frac{a}{2})$, as a function of $(a, b)$ has no critical point in the interior of the domain $\left[\frac{2}{13},\frac{1}{5}\right]\times\left[\frac{4}{2+3a},\frac{2}{1+a}\right]$. Thus, the minimum value of $f(a, b)$ is 
achieved on the boundary of the domain.
Furthermore, a series of calculations shows that $f$ is positive on the boundary. 
Consequently $f(a,b)>0$ for all $(a,b)\in\left[\frac{2}{13},\frac{1}{5}\right]\times\left[\frac{4}{2+3a},\frac{2}{1+a}\right]$. Similarly, we can show  $L(x,a,b):=I'(x)>0$ for all
%$$L(x,a,b):=I'(x)=h(x)-\left(1-x-\frac{2}{b}\right)h'(x%)-ak(x)-a\left(1+x-\frac{1}{b}-2a\right)k'(x)$$ 
$(x, a, b) \in [-\frac{1}{10},0]\times[\frac{2}{13}, \frac{1}{5}]\times\left[\frac{20}{13},\frac{26}{15}\right]$ which contains
the compact set $\left[-\frac{a}{2},0\right]\times\left[\frac{2}{13},\frac{1}{5}\right]\times\left[\frac{4}{2+3a},\frac{2}{1+a}\right]$.

For the subcase $-1+\frac{1}{b}+2a<0$ (with $a\in(\frac{2}{13}, \frac{1}{5})$) we note that  $g_{1,-2} \geq0$. In this case, a series of computations shows that  $|D |$ is given either by~\eqref{deter1} or~\eqref{deter2}. We proceed similarly to establish that the determinant is positive for $x\in [-\frac{a}{2}, 0]$.

Finally, assume $a\in[\frac{1}{5},\frac{2}{9})$. It follows that $g_{1,-2} =0$, and a series of computations shows that  $|D |$ is given  by~\eqref{deter2}, which can be shown to be positive.
Consequently,  for $(a, b)\in \Gamma'_3$, then   $|D |>0$ on $[-\frac{a}{2}, 0]$.

%$${\color{red} |D(x) |=\begin{cases}
% I(x) , & x\in [-\frac{a}{2}, 1-\frac{2}{b}+a]\\
% \frac{4a(b-1) (x+a)}{b}g_{-2,-2}(x) , & x\in [1-\frac{2}{b}+a, 0]
%\end{cases}}.$$
%In either case, we can show that $|D |>0$. 

We now consider $(a, b)\in \Gamma_3''$. Assume that $a\in (\frac{2}{7}, \frac{1}{3}]$, $b\in( \frac{4}{2+3a}, \frac{3}{2}]$. If $-1+\frac{1}{b}+a>0$, then 
\begin{equation}\label{deter5}
 |D(x)|=\begin{cases}
  J(x)  & x\in [-\frac{a}{2}, 1-\frac{1}{b}-a)\\
 g_{-2,-2}(x)   |D^{\{2,3,4\}}(x)  | & x\in [1-\frac{1}{b}-a, 0]
     \end{cases}.
\end{equation}
%In the same way we established that $I(x)>0$, we can show %that $G(x)>0.$  
As in the case of $I $, we prove that $J \neq0$ and $ |D^{\{2,3,4\}}  |>0$, since 
%\begin{equation*}
%J(x)= |\begin{array}{cccc}
% 1-x-\frac{2}{b}+2a & 1-x-\frac{2}{b}+a & 0 & 0\\
% 1+x+\frac{1}{b}-2a & 1-x-\frac{1}{b}+a & %1-x-\frac{1}{b} & 1-x-\frac{1}{b}-a\\
% 1+x-2a & 1+x-a & 1+x & 1-x-a\\
% 0 & 0 & 1+x-\frac{1}{b} & 1+x-\frac{1}{b}+a
% \end{vmatrix} \, ,
% \end{equation*} 
%\begin{equation*}
%D^{\{2,3,4\}}(x)=\begin{vmatrix}
%  1-x-\frac{1}{b}+a & 1-x-\frac{1}{b} & %1-x-\frac{1}{b}-a\\
%  1+x-a & 1+x & 1-x-a\\
 % 0 & 1+x-\frac{1}{b} & 1+x-\frac{1}{b}+a\\ 
%\end{vmatrix}>0 \, .
%\end{equation*}
$$\begin{vmatrix}
g_{-1,0}  & g_{-1,1} \\
g_{1,0}  &g_{1,1} \\ 
\end{vmatrix}>0, \quad \begin{vmatrix}
g_{-1,-1}  &g_{-1,0} \\
 0 &g_{1,0}  \\ 
\end{vmatrix}>0, \mbox{and} \begin{vmatrix}
 g_{-1,-1}  & g_{-1,1} \\
  0 &g_{1,1} \\ 
 \end{vmatrix}>\begin{vmatrix}
g_{-1,0}  & g_{-1,1} \\
g_{1,0}  &g_{1,1} \\ 
\end{vmatrix}+\begin{vmatrix}
g_{-1,-1}  &g_{-1,0} \\
 0 &g_{1,0} \\ 
\end{vmatrix}.$$ 
The same decomposition is obtained for $b\in(\frac{4}{2+3a}, \frac{3}{2}]$, $-1+\frac{1}{b}+a\geq0$ and
 $a\in (\frac{2}{7}, \frac{1}{3}]$, $b\in( \frac{3}{2};\frac{2}{1+a}]$. 
 
Let $a\in (\frac{1}{3}, \frac{2}{5}]$ and $2a-\frac{1}{b}>0$. Then,
\begin{equation}\label{deter6}
 |D(x)|=\begin{cases}
  J(x)& x\in [-\frac{a}{2},1-\frac{1}{b}-a)\\
  M(x)& x\in [1-\frac{1}{b}-a, 1-\frac{2}{b}+a)\\
 g_{-2,-2}(x)  |D^{\{2,3,4\}}(x) | & x\in [1-\frac{2}{b}+a), 0].\\
     \end{cases}
\end{equation}
%\begin{equation*}
%M(x)= |\begin{array}{cccc}
% 1-x-\frac{2}{b}+2a & 1-x-\frac{2}{b}+a & 0 & 0\\
% 1+x+\frac{1}{b}-2a & 1-x-\frac{1}{b}+a & %1-x-\frac{1}{b} & 0\\
% 1+x-2a & 1+x-a & 1+x & 1-x-a\\
% 0 & 0 & 1+x-\frac{1}{b} & 1+x-\frac{1}{b}+a. 
% \end{vmatrix} \, .
%\end{equation*}

%\begin{equation*}
%D^{\{2,3,4\}}(x)=\begin{vmatrix}
 % 1-x-\frac{1}{b}+a & 1-x-\frac{1}{b} & 0\\
 % 1+x-a & 1+x & 1-x-a\\
 % 0 & 1+x-\frac{1}{b} & 1+x-\frac{1}{b}+a
 %\end{vmatrix}>0\ \mbox{since}\
%\end{equation*}
The previous determinants are obtained in the case $a\in (\frac{1}{3}, \frac{2}{5}]$ and $2a-\frac{1}{b}\leq0$.

Let $a\in( \frac{2}{5}, \frac{1}{2})$ and $2-\frac{2}{b}-a\leq0$. Then 
\begin{equation} \label{deter7}
 |D(x)|=\begin{cases}
  g_{1,1}(x)   |D^{\{1,2,3\}}(x) |  & x\in [-\frac{a}{2}, -1+\frac{1}{b}]\\
  N(x)&x\in (-1+\frac{1}{b},1-\frac{2}{b}+a)\\
 g_{-2,-2}(x)   |D^{\{2,3,4\}}(x) | & x\in [1-\frac{2}{b}+a, 0]\\
     \end{cases}.
\end{equation}
%$$D^{\{1,2,3,4\}}(x)=\begin{vmatrix}
% 1-x-\frac{2}{b}+2a & 1-x-\frac{2}{b}+a & 0\\
% 1+x+\frac{1}{b}-2a & 1-x-\frac{1}{b}+a & %1-x-\frac{1}{b} \\
% 0 & 1+x-a & 1+x 
% \end{vmatrix} $$ and 
%$$
%N(x)=  |\begin{array}{cccc}
% 1-x-\frac{2}{b}+2a & 1-x-\frac{2}{b}+a & 0 & 0\\
% 1+x+\frac{1}{b}-2a & 1-x-\frac{1}{b}+a & %1-x-\frac{1}{b} & 0\\
% 0 & 1+x-a & 1+x & 1-x-a\\
% 0 & 0 & 1+x-\frac{1}{b} & 1+x-\frac{1}{b}+a\\ 
% \end{vmatrix}
%$$
Let $a\in( \frac{2}{5}, \frac{1}{2})$ and $2-\frac{2}{b}-a>0$. Then
\begin{equation} \label{deter8}
 |D(x)|=\begin{cases}
  P(x) & x\in [-\frac{a}{2}, 2a-1]\\
  J(x) & x\in (2a-1, 1-\frac{1}{b}-a)\\
  M(x)&x\in [1-\frac{1}{b}-a,1-\frac{2}{b}+a)\\
 g_{-2,-2}(x)   |D^{\{2,3,4\}}(x) | & x\in [1-\frac{2}{b}+a, 0]\\
     \end{cases}.
\end{equation}
%where
%$$
%P(x)=  |\begin{array}{cccc}
% 1-x-\frac{2}{b}+2a & 1-x-\frac{2}{b}+a & 0 & 0\\
% 1+x+\frac{1}{b}-2a & 1-x-\frac{1}{b}+a & %1-x-\frac{1}{b} & 1-x-\frac{1}{b}-a\\
% 0 & 1+x-a & 1+x & 1-x-a\\
% 0 & 0 & 1+x-\frac{1}{b} & 1+x-\frac{1}{b}+a\\ 
% \end{vmatrix}>0 \\ 
%$$

Similar techniques used for $ |D^{\{2,3,4\}} |$  can be use to prove that $ |D^{\{1,2,3\}} |>0$, $N>0$ and $P>0$. We conclude that for 
$(a, b)\in\Gamma'''_3$, $D$ is invertible on $ [-\frac{a}{2}, 0]$. 
\end{proof}

\begin{figure}
    \centering
    \includegraphics[width=.48\textwidth]{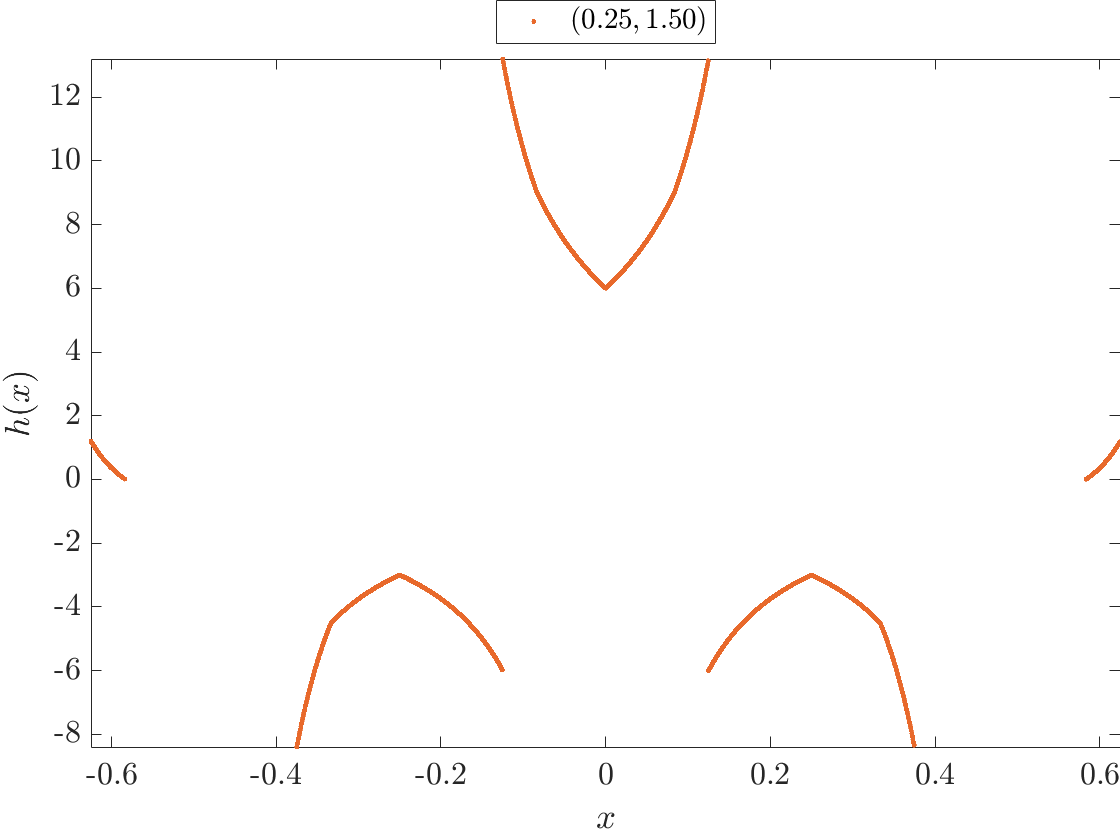}\hfill \includegraphics[width=.48\textwidth]{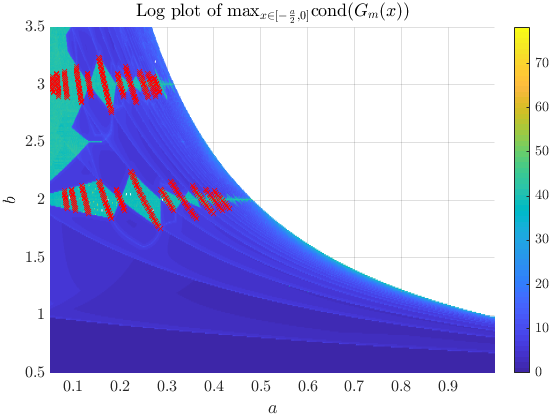}
    \caption{Although the region $\{a\in [\frac{2}{9}, \frac{2}{7}], b\in [\frac{4}{2+3a},\frac{2}{1+a}], b>1\}$ is not covered by Theorem \ref{thm-meta}, numerical evidence suggests that this set is included in $\mathcal{F}(B_2)$. On the left is a plot of nonzero values of the dual $h$ of $\mathcal{G}(B_2,a,b)$ corresponding to $a=0.25,b= 1.50$, and on the right is a plot of the condition number of the matrices $G_m$, indicating the invertibility of the matrices in this region. The red marks indicate known points $(a,b)\not \in \mathcal{F}(B_2)$.}
    \label{fig:gap}
\end{figure}

\begin{rem}
Using numerical simulations (Figure \ref{fig:gap}) we observe  that for  $(a, b) \in \Gamma'''_3:=\{(a,b)\in\mathbb{R}^2_+:a\in (\frac{2}{9}, \frac{2}{7}], b\in (\frac{4}{2+3a})\}$, then $|D|>0$. However, we have not yet been able to prove this.
\end{rem}

%\begin{proof} 
%Let $\{D_k(x)\}_{k=1}^4$ denote the columns of $D(x)$. We claim %that $\{D_k\}_{k=1}\subset L^2(\R, \C^4)$ is a linearly %independent set where $L^2(\R, \C^4)$ is a vector-valued Hilbert %space with the inner product $$\ip{f}{g}_{L^2(\R, \C^4)}=\int_\R %\sum_{k=1}^4 f_k(x)\overline{g_k(x)}dx$$ for $f=(f_k)_{k=1}^4, %g=(g_k)_{k=1}^4 \in L^2(\R, \C^4)$. To establish this, it %suffices to show that the Gram matrix $M=\begin{pmatrix} %\ip{D_k}{D_\ell}_{L^2(\R, \C^4)}\end{pmatrix}_{k, \ell=1}^4$ is %positive definite. A series of computations shows that %$$M=4\begin{pmatrix} \hat{h}(0) & \hat{h}(a) & \hat{h}(2a) & %\hat{h}(3a)\\
%\overline{\hat{h}(a)} &\hat{h}(0) & \hat{h}(a) & \hat{h}(2a)\\
%\overline{\hat{h}(2a)} & \overline{\hat{h}(a)}& \hat{h}(0) & %\hat{h}(a)\\
%\overline{\hat{h}(3a)} & \overline{\hat{h}(2a)}& %\overline{\hat{h}(a)}& \hat{h}(0)
%\end{pmatrix}$$
%where $h(x)=|\hat{B_{2}}(x)|^2 \in L^1(\R)$. Since $\hat{h}$ is %a positive definite function, we conclude that $M$ is a positive %definite matrix which establish our claim that the columns of %$D(x)$ are linearly independent. hence, 

%\end{proof}

\begin{proof}[Proof of Theorem~\ref{thm2:b2new}]
We recall from Propositions~\ref{prop:invgamma_3} that for $(a, b)\in\Gamma_3$, we have $|G_3|\neq 0$ on $[-\frac{a}{2}, 0]$.  Since $h$ is compactly supported, we know that $h(x)=0$ for $|x|>\frac{5a}{2}$. In addition, on  $\left[-\frac{5a}{2},\frac{5a}{2}\right]$, $h$ is determined by the entries of   $H_3=b(G^{-1}_3)_3$, where $H_3$ is given by \eqref{eq:H}.
% $$ 
% \begin{pmatrix}
%  h(x-2a) \\ h(x-a)\\ h(x)\\ h(x+a) \\ h(x+2a) \end{pmatrix}=G^{-1}_3(x)\begin{pmatrix}0\\ 0\\ b\\0 \\0\end{pmatrix}=b(G^{-1}_3(x))_3$$ 

Then, we can determine $h$ on $(-\frac{a}{2}+ka, ka]$ where \mbox{$k=0, \pm 1, \pm 2$}. This, along with the symmetry of $h$, will define the function $h$ everywhere except possibly at $\frac{ka}{2}$ for $k=\pm1, \pm2, \pm3, \pm 5.$   In particular, $h(x+2a)=0$ for $x\in (-\frac{a}{2}, 0)$. That is, $h=0$ on $(-\frac{5a}{2}, -2a)$, and $(2a, \frac{5a}{2})$.

 Consequently,  $h$  is a compactly supported and bounded function for which $\cG(h,a,b)$ is a Bessel sequence. By construction,  it also follows that  $g$ and $h$ are dual windows. 
\end{proof}

\begin{figure}
    \centering
\includegraphics[width=.5\linewidth]{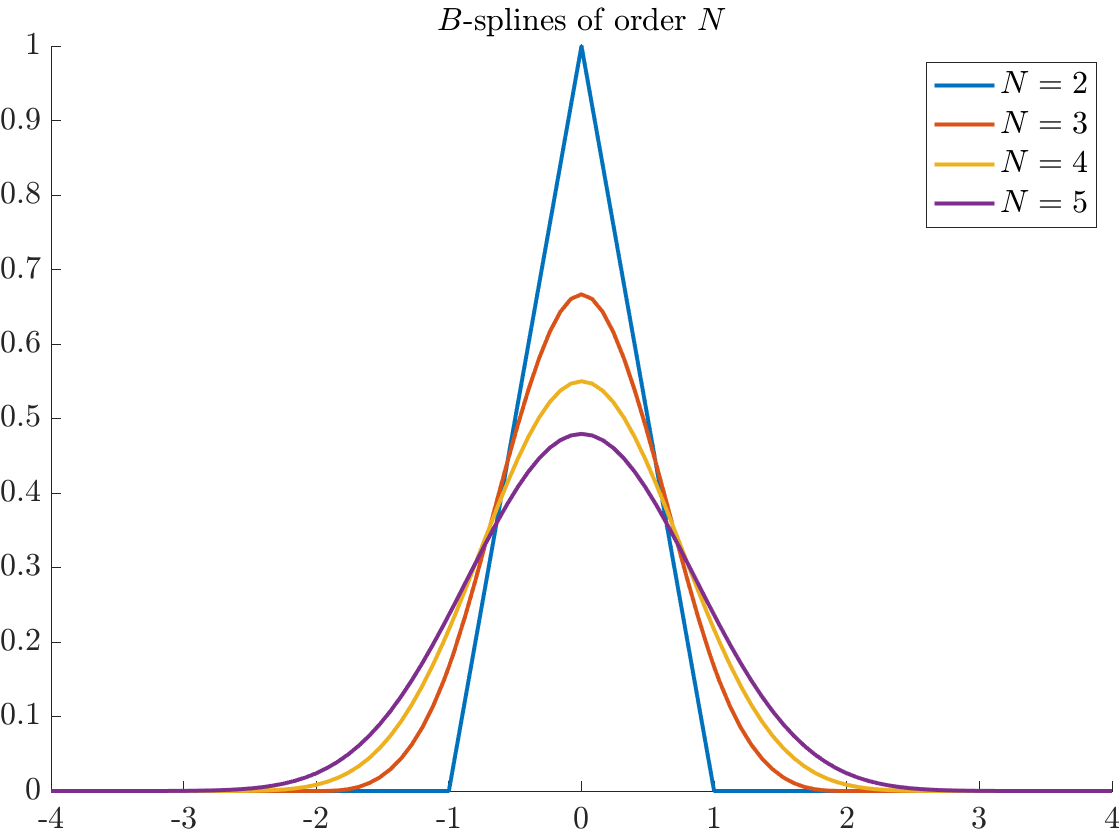}\includegraphics[width=.5\linewidth]{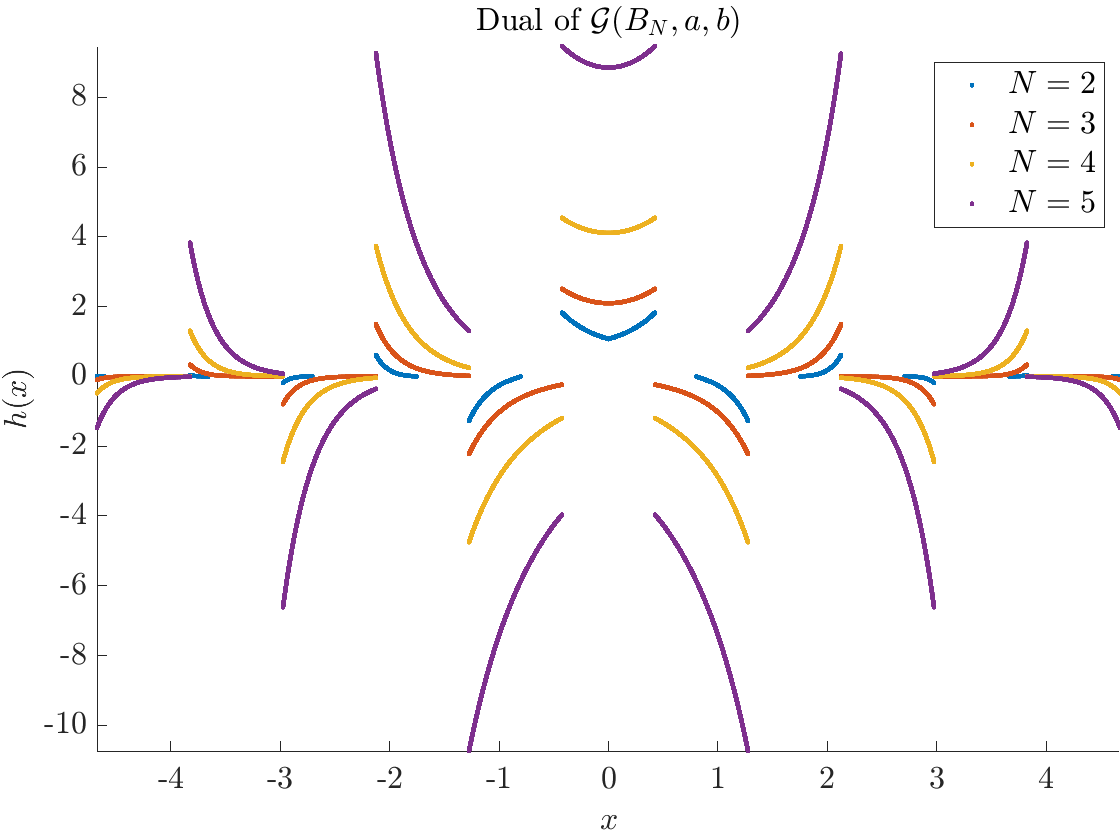}
    \caption{ The results presented here extend in a straightforward way to higher-order B-splines (left). On the right are plots of the nonzero values of the corresponding dual frames $h$ of $\mathcal{G}(B_N, a, b)$ for $(a,b)= (0.85, 1.05)$ and $N=2,3,4,5$.}
    \label{fig:BN}
\end{figure}

\begin{rem}
\begin{enumerate}
    \item We observe that the dual window $h$ constructed from Theorem~\ref{thm1:b2new} and Theorem~\ref{thm2:b2new} is discontinuous. This is proved in the same way as in \cite[Remark 4]{AtiKouOko1}.
    \item On the line $b=\frac{3}{2}$, all of the dyadic points in the set  $\{(\frac{1}{2^j}, \frac{3}{2})\}_{j=1}^\infty \cup \{(\frac{3}{2^j}, \frac{3}{2})\}_{j=3}^\infty$ belong to  $\mathcal{F}(B_2)$, since most of these points belong $\Gamma_3$, except  for $(\frac{1}{4},\frac{3}{2}), (\frac{3}{8}, \frac{3}{2}),$ and $ (\frac{1}{2},\frac{3}{2})$ for which we omit the proof. 
    \item Our results extend to regions beyond the frame set established here and also to higher-order B-splines, as demonstrated in Figures \ref{fig:gap} and \ref{fig:BN}.
\end{enumerate}

\end{rem}

\section*{Acknowledgment} Part of this work was completed while the first-named author was a visiting graduate student in  the Department of Mathematics at the University of Maryland during the Fall 2017 semester. He would like to thank the Department for its hospitality and the African Center of Excellence in Mathematics and Application (CEA-SMA) at the  Institut de Math\'ematiques et de Sciences Physiques (IMSP) for funding his visit. K.~A.~Okoudjou  was partially supported by a grant from the Simons Foundation $\# 319197$, by ARO grant W911NF1610008, and the National Science Foundation grant DMS 1814253.

%%%%%%%%%%%%%%%%%%%%%%%

\bibliographystyle{amsplain}
\bibliography{B2FS}
\newpage
\section{Appendix}
\tiny
\begin{lstlisting}[caption = {MATLAB code for plotting the dual of $\mathcal{G}(B_N, a,b).$},style=Matlab-editor,
basicstyle=\ttfamily,
escapechar=]
function plotDual
a = .85;
b = 1.05;
N = 2;
g = @(x) fnval(spmak(-N:N,[0 1 0]),x);
x=-a/2:.001:a/2;

% Set up and solve the linear systems for h(x)
for j=1:length(x)
    [m, Gm]=G(x(j),g, a, b);
    if m==0
        H = NaN;
    else
        bv = zeros(2*m-1,1); 
        bv(m) = b;
        H(:,j) = Gm\bv;
    end
end

% concatenate the components of H
X = []; h = [];
for k = -(m-1):(m-1)
    X = [X x+k*a];
    h = [h, H(k+m,:)];
end

figure;
plot(X(abs(h)>1e-10), h(abs(h)>1e-10), '.')
title(sprintf('Dual of $\\mathcal{G}(B_N)$ \n N=%d, a=%1.2f, b=%1.2f,',N, a, b))
end

function [m, Y] = G(x, g, a, b, m_max)
% Returns the (2m-1) x (2m-1) matrix Y = G_3(x) and the integer m between 1 and m_max such that (a,b) is in \Lambda_m.
% If (a,b) is not in \Lambda_m for any m\geq 1, Y=NaN.
% The window g is a function handle, m_max is a large integer and x, a, and b are real numbers.

if nargin<5
    m_max=50;
end

m = NaN;
for mm = 1:m_max
    if (b>2*(mm-1)/(2+(2*mm-3)*a)) && (b<2*mm/(2+(2*mm-1)*a)) && (b<2/(1+a)) && (b>1)
         m= mm;
    end
end
if ~isnan(m)
    l = meshgrid(-(m-1):(m-1));
    Y = g(x-l'/b+l*a);
else
    Y=NaN;
    fprintf('Error. Perhaps try again with larger m_max.')
end
end

\end{lstlisting}
\end{document}